\documentclass[12pt,reqno]{amsart}
\usepackage{graphics}
\usepackage{tikz}
\usetikzlibrary{positioning,decorations.pathreplacing}
\usepackage{graphicx} 
\usepackage{psfrag} 
\usepackage{amssymb,a4wide,amsthm,amsmath,amsfonts}
\usepackage{mathrsfs}
\usepackage{color}

\theoremstyle{plain}
\newtheorem{theorem}{Theorem}[section]
\newtheorem{fact}[theorem]{Fact}
\theoremstyle{remark}

\theoremstyle{plain}

\numberwithin{equation}{section}

\newcommand{\R}{\mathbb{R}}
\newcommand{\N}{\mathbb{N}}
\newcommand{\Z}{\mathbb{Z}}

\newcommand{\F}{\mathcal{F}}

\newtheorem*{main}{Theorem (Main result)}

\newcounter{probcount}

\begin{document}
\title{On Lipschitz-free spaces over spheres of Banach spaces}
\author{Pedro L. Kaufmann}
\address{Instituto de Ci\^encia e Tecnologia da Universidade federal de S\~ao Paulo, Av. Cesare Giulio Lattes, 1201, ZIP 12247-014 S\~ao Jos\'e dos Campos/SP, Brasil}
\thanks{Grants 2016/25574-8 and 2017/18623-5, S\~ao Paulo Research Foundation (FAPESP)} 
\author{Leandro Candido}
\begin{abstract}
We prove that, for each Banach space $X$ which is isomorphic to its hyperplanes, the Lipschitz-free spaces over $X$ and over its sphere are isomorphic.
\end{abstract}
\maketitle

\section{Introduction}

Geometry of Lipschitz-free spaces have become a very active research subject in recent years. An approach to the topic that has been revealed as very fruitful consists in establishing relations between the free space over a metric space and the free spaces over certain of its subsets. Let us mention some examples: \medskip

$\bullet$ Godefroy and Ozawa showed in \cite{godefroy2014free} that, given a separable Banach space $X$ and a convex closed subset $C$, the free space $\F(C)$ contains an isometric and 1-complemented copy of $\F(X)$. With this, they were able to construct a free space over a compact metric space failing to have the approximation property. \medskip

$\bullet$ In \cite{Kalton}, Kalton shows that the free space over a pointed metric space $(M,d,0)$ admits a complemented copy in the $\ell_1$-sum of $\F(M_k),k\in\Z$, where $M_k$ are the annuli defined by $M_k=\{x\in M: 2^{k-1}\leq d(x,0)\leq 2^{k+1}\}$. This result has a range of applications; in the same paper it is used for instance to prove that the free space over uniformly discrete metric spaces has the approximation property. Another application is pointed out in the next item. \medskip

$\bullet$ In \cite{kaufmann2015products} it is shown that the free space over a Banach space is always isomorphic to the free space over its ball, using Kalton's aforementioned result. This result was for instance used in its finite dimensional version by C\'uth, Doucha and Wojtaszczyk in \cite{cuth2016structure} to prove that $\F(M)$ is weakly sequentially complete for every subset $M$ of $\R^n$. \medskip

$\bullet$ More recently, Aliaga, No\^us, Petitjean and Proch\'azka obtained a related result that can be found in \cite{aliaga2020compact}: the free space over a complete metric space $M$ is weakly sequentially complete if and only if $\F(K)$ is weakly sequentially complete, for each compact subset $K$ of $M$. \medskip

Motivated especially by the second and third examples, one can ask: \medskip

\textbf{Question.} Given a Banach space $X$, how can we relate the free space over it and over its sphere, $S_X$? \medskip

As we shall see, using well known facts we can show that we can always find a complemented copy of $\F(S_X)$ in $\F(X)$. When we assume that, moreover, $X$ is isomorphic to one of its hyperplanes (and therefore to all of its hyperplanes \cite[Exercice 2.9]{fabian2011banach}), we have more: 

\begin{main}\label{main}
Let $X$ be a Banach space such that $X$ is isomorphic to its hyperplanes. Then $\F(X)$ and $\F(S_X)$ are isomorphic.
\end{main}
 
 In the remainder of this text we will present a short proof of this  statement. 
 
 \section{Proof of main result}
 
 Throughout the text, we shall adopt the following notation. Given $X$ and $Y$ Banach spaces, $X\simeq Y$ stands for \emph{$X$ is isomorphic to $Y$}. $X\cong Y$ stands for \emph{$X$ is isometrically isomorphic to $Y$}. $X\stackrel{c}{\hookrightarrow} Y$ stands for \emph{$X$ is isomorphic to a complemented subspace of $Y$}. The closed unit ball of $X$ is denoted by $B_X$.

 Let us start recall some basic properties and notions orbiting around free spaces which we will use. For the reader unfamiliar to the subject, we refer to Weaver's book \cite{weaver1999lipschitz} (where Lipschitz-free spaces are denominated Arens-Eells spaces), Godefroy and Kalton's paper \cite{godefroy2003lipschitz} and Godefroy's survey \cite{godefroy2016survey} for a thorough introduction. 

Given a pointed metric space $(M,0)$, we denote by $\mathrm{Lip}_0(M)$ the Banach space of real valued Lipschitz functions defined on $M$ vanishing at 0. For each $x\in M$, we define an evaluation functional $\delta(x)\in \mathrm{Lip}_0(M)^*$ by setting $\delta(x)(f)=f(x)$.  The Lipschitz-free space over $M$, denoted by $\F(M)$, is defined as the closed linear span of all evaluation functionals in $\mathrm{Lip}_0(M)^*$. $\delta_M = \delta: x\in M\mapsto\delta(x)\in \F(M)$ is easily seen to be an isometric embedding.  $\F(M)$ is a canonical isometric predual for $\mathrm{Lip}_0(M)$, giving it a $w^*$ topology which coincide with the topology of pointwise convergence  on bounded subsets of $\mathrm{Lip}_0(M)$.

A powerful property of this spaces is that they allow linear interpretation of Lipschitz maps between metric spaces in a unique fashion. Namely, given two pointed metric spaces $(M,0_M)$, $(N,0_N)$ and a Lipschitz map $f:M\to N$ taking $0_M$ to $0_N$, there exists a unique bounded linear operator $T:\F(M)\to \F(N)$ such that $T\circ \delta_M=\delta_N\circ f$. The norm of the operator $T$ coincides with the Lipschitz constant of $f$. In particular, if $M$ and $N$ are bi-Lipschitz equivalent (that is, if there exists a Lipschitz bijection with Lipschitz inverse between them), then $\F(M)\simeq \F(N)$. 

Now suppose that $0\in N\subset M$. Whenever we can find a bounded linear operator $E: \mathrm{Lip}_0(N) \to \mathrm{Lip}_0(M)$ such that $E(f)|_{N} = f$ (that is, $E$ is an \emph{extension} operator) which is moreover $w^*$ continuous, it admits a preadjoint $E_*: \F(M)\to \F(N)$ which can actually be realized as a (bounded linear) projection if we see $\F(N)$ as a subspace of $F(M)$. It follows that $\F(N)\stackrel{c}{\hookrightarrow}\F(M)$.
There are several ways to pursue such extension operators (the interested reader should refer to \cite{BruBru} and \cite{LN05}). One of the simplest (and which shall suffice for our purposes) is to look for a Lipschitz projection from $M$ onto $N$. If such projection $P:M\to N$ exists, we can define $E$ by putting $E(f)=f\circ P$.  \medskip

We recall at this point that  Benyamini and Sternfeld proved in \cite{benyamini1983spheres} that there exists a Lipschitz retraction from $B_X$ onto $S_X$, whenever $X$ is an infinite dimensional Banach space. It is straightforward to obtain a retraction from $X$ onto $S_X$ by composing it with the Lipschitz projection $G: X\to B_X$ given by the formula \[G(x)= \left\{
\begin{array}{ll}
\frac{x}{\|x\|} & \text{ if }\|x\|\geq 1;\\
x & \text{ if } \|x\|\leq 1.
\end{array} \right.\]
By the above discussion we can conclude the the following holds: 

\begin{fact}\label{Fact2} Let $X$ be an infinite dimensional Banach space. Then,  $\F(S_X)\stackrel{c}{\hookrightarrow}\F(X)$. 
\end{fact}

The following observation will also be useful: 

\begin{fact}\label{Fact1} If $X$ and $X$ are isomorphic Banach spaces, then their spheres are bi-Lipschitz equivalent. In particular, $\F(S_X)\simeq \F(S_Y)$. 

\end{fact}
 Indeed, suppose that $T:X\to Y$ is an isomorphism. Then the function $f:S_X\to S_Y$ given by the formula $f(x)=T(x)/\|T(x)\|$ is a bi-lipschitz isomorphism with distortion less or equal to $(2\|T^{-1}\|\|T\|)^2$. We are now in position to prove our main result:

\begin{proof}[Proof of main result] 
Since $X$ is isomorphic to its hyperplanes, $X$ is isomorphic to $\R\oplus_\infty X=Y$. In particular, by Fact \ref{Fact1}, $S_X$ and $S_Y$ are Lipschitz equivalent. Note that $D=\{(1,x):x\in B_X\}$ is a subset of $S_Y$, isometrically isomorphic to $B_X$. Define $\pi: S_Y\to D$ by $\pi(\lambda,x)=(1,x)$. $\pi$ is a 1-Lipschitz retract. Indeed, it clearly fixes points of $D$, and for each $(\lambda,x),(\alpha,y)\in S_Y$,
\begin{align*}
\|\pi(\lambda,x)-\pi(\alpha,y)\|_Y&=\|(0,x-y)\|_Y=\|x-y\|_X\\
&\leq \max\{|\lambda-\alpha|,\|x-y\|_X\} = \|(\lambda,x)-(\alpha,y)\|_Y.
\end{align*}
It follows that 
\begin{eqnarray}
\F(D)\stackrel{c}{\hookrightarrow}\F(S_Y).
\label{FD}
\end{eqnarray}
 Recall that, by \cite[Corollary 3.3]{kaufmann2015products}, $\F(D)\cong \F(B_X)\simeq \F(X)$. By Fact \ref{Fact1}, $\F(S_X)\simeq \F(S_Y)$. We then may rewrite (\ref{FD}) as $\F(X)\stackrel{c}{\hookrightarrow}\F(S_X)$.

On the other hand, by Fact \ref{Fact2}, $\F(S_X)\stackrel{c}{\hookrightarrow}\F(X)$. An application of Pe\l czy\'nski decomposition method leads us to the conclusion.
\end{proof}

We remark that recently Albiac, Ansorena, C\'uth and Doucha proved that $\F(\R^n)\simeq \F(S_{\R^{n+1}})$, for each $n\in \N$ (see \cite{albiac2020lipschitz}). Of course, $\R^n$ does not satisfy the criteria to apply our main result. In view of our discussion, this relates to the question of whether or not $\F(\R^n)\not\simeq \F(\R^m)$ when $n\neq m$, which is open for $n,m\geq 2$, as far as we know. For $n=1$, it is well known that $\F(\R^2)\not \hookrightarrow \F(\R)$ due to a result of Naor and Schechtman (\cite{naor2007planar}).

\bibliographystyle{siam}
\bibliography{references}
\end{document}